\documentclass{article}
\catcode`\@=11
\newtheorem{theorem}{Theorem}
\newtheorem{lemma}{Lemma}

\newtheorem{example}{Example}
\newtheorem{remark}{Remark}

\def\demo{\noindent{\bf Proof .-}}
\def\section{\@startsection {section}{1}{\z@}{-3.5ex plus -1ex
minus-.2ex}{2.3ex plus .2ex}{\normalsize\bf}}

\def\bz{\hbox{\it Z\hskip -4pt Z}}

\newcommand{\het}{H_{\rm et}}
\newcommand{\hc}{H_{\rm c}}

\begin{document}
\begin{center}
{\Large\bf \textsc{Certain minimal varieties are set-theoretic complete intersections}}
\end{center}
\vskip.5truecm
\begin{center}
{Margherita Barile\footnote{Partially supported by the Italian Ministry of Education, University and Research.}\\ Dipartimento di Matematica, Universit\`{a} di Bari,Via E. Orabona 4,\\70125 Bari, Italy}
\end{center}
\vskip1truecm
\noindent
{\bf Abstract} We present a class of  homogeneous ideals which are generated by monomials and binomials of degree two and are set-theoretic complete intersections. This class includes certain reducible varieties of minimal degree and, in particular,  the presentation ideals of the fiber cone algebras of monomial varieties of codimension two. 

\section*{Introduction and Preliminaries}
Let $K$ be an algebraically closed field, and let $R$ be a polynomial ring  in $N$ indeterminates over $K$. Let $I$ be a proper reduced ideal of $R$ and consider the variety $V(I)$ defined in the affine space $K^N$ (or in the projective space ${\bf P}_K^{N-1}$, if $I$ is homogeneous and different from the maximal irrelevant ideal) by the vanishing of all polynomials in $I$. By Hilbert Basissatz there are  finitely many polynomials $F_1,\dots, F_s\in R$ such that $V(I)$ is defined by the equations $F_1=\cdots= F_s=0$. By Hilbert Nullstellensatz this is equivalent to the ideal-theoretic condition

\begin{equation}\label{1}I=\sqrt{(F_1,\dots, F_s)}.\end{equation}
\noindent
Suppose $s$ is minimal with respect to this property. It is well known that height\,$I\leq s$. If equality holds, $I$ is called a {\it set-theoretic complete intersection} (s.t.c.i.) on $F_1,\dots, F_s$. \newline
Exhibiting significant examples of s.t.c.i. ideals (or, more generally, determining the minimum number of equations defining given varieties, the so-called {\it arithmetical rank}) is one of the most difficult problems in algebraic geometry. The main problem is finding polynomials $F_1,\dots, F_s$ fulfilling (\ref{1}), which need not be part of a minimal generating system for 
$I$. This task cannot be accomplished by a general constructive method. There are, however, a few results that allow us to settle several special cases. One of these is due to Schmitt and Vogel.
\begin{lemma}\label{Schmitt}{\rm [\cite{SV}, p.~249]} 
Let $P$ be a finite subset of elements of $R$. Let $P_0,\dots, P_r$ be subsets of $P$ such that
\begin{list}{}{}
\item[(i)] $\bigcup_{l=0}^rP_l=P$;
\item[(ii)] $P_0$ has exactly one element;
\item[(iii)] if $p$ and $p''$ are different elements of $P_l$ $(0<l\leq r)$ there is an integer $l'$ with $0\leq l'<l$ and an element $p'\in P_{l'}$ such that $(pp'')^m\in(p')$ for some positive integer $m$.
\end{list}
\noindent
We set $q_l=\sum_{p\in P_l}p^{e(p)}$, where $e(p)\geq1$ are arbitrary integers. We will write $(P)$ for the ideal of $R$ generated by the elements of $P$.  Then we get
$$\sqrt{(P)}=\sqrt{(q_0,\dots,q_r)}.$$
\end{lemma}
This result, together with its refinements and generalizations established in \cite{B2}, is especially useful for ideals generated by monomials. An interesting class of s.t.c.i monomial ideals was introduced by Lyubeznik \cite{L} and studied in \cite{B1}. 
\newline
Big classes of s.t.c.i. ideals generated by binomials have been characterized among the toric ideals (see, e.g., \cite{BL}, \cite{BMT}, \cite{K} and \cite{Br}, \cite{El}, \cite{E}, \cite{P} for toric curves), and some cases are also known among the determinantal ideals  of two-row matrices. One of these was treated by Robbiano and Valla \cite{RV1}, another one by Bardelli and Verdi \cite{BV}.  In this paper we generalize Bardelli and Verdi's result by presenting a class of s.t.c.i. ideals generated by some minors  and some products of entries of certain blockwise defined matrices, which were considered in \cite{BM1}. The corresponding varieties include the fiber cones of codimension two monomial varieties, which were studied by Gim\'enez, Morales and Simis in \cite{G} and  \cite{GMS}, and, furthermore, they belong to the (larger) class of reducible varieties of minimal degree classified by Xamb\'o \cite{X}. These are all defined by monomials and binomials, but they are not all s.t.c.i.. At the end we shall exhibit a counterexample, which is also interesting from another point of view: in positive characteristics its arithmetical rank is strictly greater than its cohomological dimension, which seems to be a rare property.  
\section{A class of set-theoretic complete intersections}
Let $r$ be a positive integer and consider the two-row matrix
\begin{equation}\label{2}A=\left(B_1\|B_2\|\dots\|B_r\right),\end{equation}
\noindent
where, for all $i=1,\dots, r$, $B_i$ is the $2\times c_i$-matrix
\begin{equation}\label{2'}B_i=\left(\begin{array}{ccccc}
X^i_1&X^i_2&\dots&X^i_{c_i-1}&X^i_{c_i}\\
X^i_2&X^i_3&\dots&X^i_{c_i}&X^i_{c_i+1}
\end{array}
\right).
\end{equation}
\noindent
Here \underline{$X$}$=\{X^i_j\}$ is a set of $N$ indeterminates over $K$, and $X^i_j\ne X^h_k$ for $(i,j)\ne(h,k)$, with one only possible exception: for every index $i$, $1\leq i<r$, there is at most one index $i'>i$ such that $X^i_{c_i+1}=X^{i'}_1$. This kind of matrix was introduced by Gim\'enez \cite{G}, and also considered in \cite{BM1}, where $A$ was called a {\it barred matrix} and the $B_i$'s were called the {\it big blocks} of $A$. 
\begin{example}\label{Example1}{\rm  An example of matrix of type (\ref{2}) in the indeterminates $X_1,\dots, X_{11}$ is the following:
$$A=\left(\begin{array}{ccccccccccccc}
X_1&\vline\, \vline&X_3&X_4&\vline\, \vline&X_5&X_6&X_7&\vline\, \vline&X_9&X_{10}\\
X_2&\vline\, \vline&X_4&X_5&\vline\, \vline&X_6&X_7&X_8&\vline\, \vline&X_{10}&X_{11}
\end{array}\right)
$$}
\end{example}
Fix one index $i$, $1\leq i\leq r$. If $c_i>1$, for all indices $j,j'$, $1\leq j<j'\leq c_i$ let $M^i_{jj'}$ be the minor of $B_i$ formed by the $j$-th and the $j'$-th column, i.e., 
$$M^i_{jj'}=\left|\begin{array}{cc}
X^i_j&X^i_{j'}\\
X^i_{j+1}&X^i_{j'+1}
\end{array}
\right|
=X^i_jX^i_{j'+1}-X^i_{j'}X^i_{j+1}.$$
\noindent Let $I_i\in K[X^i_1,\dots, X^i_{c_i+1}]$ be the ideal generated by all minors $M^i_{jj'}$, and consider 
\begin{equation}\label{3}
\!\!\!\!\!\!\!\! F^i_j=\displaystyle\sum_{k=0}^{j}(-1)^k{{j}\choose k}(X^i_{j+2})^{j-k}X^i_{k+1}(X^i_{j+1})^k,\quad\mbox{ for } j=1,\dots, c_i-1,
\end{equation}
\noindent
The first part of the next theorem can be found in \cite{EG}, pp.~118ff, the second part is due to Bardelli and Verdi \cite{BV} (see  \cite{RV2}, Section 2). In order to simplify the notation, we shall omit the index $i$ in the claim.
\begin{theorem}\label{theorem1} The ideal $I$ is prime of height $c-1$. It is a s.t.c.i. on $F_1, \dots, F_{c-1}$. 
\end{theorem}
The ideal $(I_i)$  thus defines an irreducible curve in ${\bf P}_K^{c_i}$, which is known as a {\it rational normal curve}, and is a special kind of {\it (rational normal) scroll}: we refer to Harris \cite{H} for an introductory treatment of this notion. A matrix of type (\ref{2'}) will be called a {\it scroll matrix}. As a consequence of Theorem \ref{theorem1} and Hilbert Nullstellensatz we also have that
\begin{equation}\label{4}\sqrt{(F^i_1,\dots,F^i_{c_i-1})}=I_i.\end{equation}
\noindent
 Following \cite{BM1}, we can associate with matrix $A$ an ideal $J$ of $R=K[$\underline{$X$}$]$. This ideal $J$ is generated by the union of
\begin{list}{}{}
\item[(i)] the set of all minors $M^i_{jj'}$  with $i=1,\dots, r$, and $1\leq j<j'\leq c_i$ (we set $M^i_{jj'}$ equal to the empty set whenever $c_i=1$);
\item[(ii)] the set of all products $X^i_jX^{i'}_{j'}$, with $1\leq i<i'\leq r$, $1\leq j\leq c_i$, $2\leq j'\leq c_{i'}+1$. These are the products of one entry of the upper row of $B_i$ and one entry of the lower row of one of the blocks $B_{i'}$ following on the right.  
\end{list}
Note that 
\begin{equation}\label{5}I_i\subset J\qquad\qquad\mbox{for all }i=1,\dots, r\end{equation}
\noindent
\begin{example}\label{Example2}{\rm The ideal $J$ associated with the matrix $A$ given in Example \ref{Example1} is generated by the following 28 elements
\begin{list}{}{}
\item[(i)] $M^2_{12}=X_3X_5-X_4^2$,\vskip.1truecm\noindent\hphantom{(i)}$M^3_{12}=X_5X_7-X_6^2$,  $M^3_{23}=X_6X_8-X_7^2$,  $M^3_{13}=X_5X_8-X_6X_7$,
\vskip.1truecm\noindent\hphantom{(i)}$M^4_{12}=X_{9}X_{11}-X_{10}^2$;
\item[(ii)]$X_1X_4$, $X_1X_5$, $X_1X_6$, $X_1X_7$, $X_1X_8$, $X_1X_{10}$, $X_1X_{11}$,
\vskip.1truecm\noindent\hphantom{(ii)}$X_3X_6$, $X_3X_7$, $X_3X_8$, $X_3X_{10}$, $X_3X_{11}$,
\vskip.1truecm\noindent\hphantom{(ii)}$X_4X_6$, $X_4X_7$, $X_4X_8$, $X_4X_{10}$, $X_4X_{11}$,
\vskip.1truecm\noindent\hphantom{(ii)}$X_5X_{10}$, $X_5X_{11}$, $X_6X_{10}$, $X_6X_{11}$, $X_7X_{10}$, $X_7X_{11}$.
\end{list}}
\end{example}
We introduce one more piece of notation. We set
\begin{equation}\label{6}G_k=\sum_{i=1}^{r-k}X^i_1X^{i+k}_{c_{i+k}+1},\qquad\mbox{for all }k=1,\dots, r-1.\end{equation}
\noindent
In other words, $G_k$ is the sum of products of the entry lying in the left upper corner of a block $B_i$ and the entry lying in the right lower corner of the block $B_{i+k}$. In particular every summand of $G_k$ is an element of the set defined in (ii), therefore
\begin{equation}\label{7} G_k\in J\qquad\mbox{for all }k=1,\dots, r-1.\end{equation}
\noindent
We are now ready to state our main result, which generalizes Theorem \ref{theorem1}.
\begin{theorem}\label{theorem2} The ideal $J$ is reduced of pure height $\sum_{i=1}^r c_i-1$. It is a s.t.c.i. on 
$$F^1_1,\dots, F^1_{c_1-1},\dots, F^r_1,\dots,F^r_{c_r-1},G_1,\dots, G_{r-1}.$$
\end{theorem}
\demo The first part of the claim is Corollary 1.3 in \cite{BM1}. We prove the second part. First of all we remark that the number of polynomials $F^i_j, G_k$ is
$$\sum_{i=1}^r (c_i-1)+r-1 = \sum_{i=1}^r c_i -r+r-1=\sum_{i=1}^r c_i-1.$$
\noindent 
Next we have to show that, in $K^n$, the variety $V$ defined by $J$ coincides with the variety $V'$ defined by the polynomials $F^i_j, G_k$. According to (\ref{4}), (\ref{5}) and (\ref{7}), all these polynomials belong to $J$, so that we certainly have $V\subset V'$. For the opposite inclusion, consider a point ${\bf x}\in K^n$ such that 
\begin{eqnarray} F^i_j({\bf x})&=&0\qquad\mbox{for all }i=1,\dots, r, j=1,\dots, c_{i}-1,\label{8}\\
G_k({\bf x})&=&0\qquad\mbox{for all }k=1,\dots, r-1.\label{9}
\end{eqnarray}
\noindent
Our claim is that all polynomials in the sets defined in (i) and (ii) vanish at ${\bf x}$.  By virtue of Theorem \ref{theorem1}, (\ref{8}) implies that, for all $i=1,\dots, r$ and $1\leq j<j'\leq c_i$,
\begin{equation}\label{10} M^i_{jj'}({\bf x})=0.\end{equation}
\noindent
 Thus there remains to prove that, for all $i,j, i', j'$ with $1\leq i<i'\leq r, 1\leq j\leq c_i, 2\leq j'\leq c_{i'}+1$,
\begin{equation}\label{11}x^i_jx^{i'}_{j'}=0,\end{equation}
\noindent
where $x^i_j$ is the monomial $X^i_j$ evaluated at ${\bf x}$.
Note that, according to (\ref{6}), 
$$G_{r-1}=X^1_1X^r_{c_r+1}.$$
\noindent
Let $1\leq k<r-1$, and consider the product of two arbitrary summands of $G_k$:
$$\Pi=X^i_1X^{i+k}_{c_{i+k}+1}\cdot X^{i'}_1X^{i'+k}_{c_{i'+k}+1},$$
\noindent
where $1\leq i<i'\leq r-k$. Let $k'$ be such that $i'+k=i+k'$. Then $k'>k$, and $X^{i}_1X^{i'+k}_{c_{i'+k}+1}=X^{i}_1X^{i+k'}_{c_{i+k'}+1}$ is a summand of $G_{k'}$ which divides $\Pi$. Let for all $k=0,\dots, r-2$, $P_k$ be the set of all summands of $G_{r-1-k}$: then Lemma \ref{Schmitt} applies. Thus (\ref{9}), together with Hilbert Nullstellensatz, implies that \begin{equation}\label{12}x^i_1x^{i+k}_{c_{i+k}+1}=0, \qquad\mbox{for all }i=1,\dots, r-1, k=1,\dots, r-i.\end{equation}
\noindent
Now suppose for a contradiction that (\ref{11}) is false, i.e., that
 for some $i,j,j',k$ with $1\leq i\leq r-1, 1\leq k\leq r-i, 1\leq j\leq c_i, 2\leq j'\leq c_{i+k}+1$,  \begin{equation}\label{13}x^i_jx^{i+k}_{j'}\neq0.\end{equation}
\noindent
For two fixed indices $i$ and $k$, let $j'-j$ be maximal with respect to (\ref{13}). Then, by (\ref{12}), necessarily $j'-j\ne c_{i+k}$, which means $j'-j<c_{i+k}$, i.e., $1<j\leq c_i$ or 
$2\leq j'<c_{i+k}+1$. First suppose that $1<j\leq c_i$. Then, by the maximality condition, we have
\begin{equation}\label{14} x^i_{j-1}x^{i+k}_{j'}=0, \end{equation}
\noindent which, in view of (\ref{13}), implies that
\begin{equation}\label{15} x^i_{j-1}=0. \end{equation}
\noindent
On the other hand, by (\ref{10}),
\begin{equation}\label{16} 
M^i_{j-1j}({\bf x})=\left|\begin{array}{cc}
x^i_{j-1}&x^i_j\\
x^i_j&x^i_{j+1}
\end{array}
\right|
=x^i_{j-1}x^i_{j+1}-({x^i_j})^2=0.
\end{equation}
\noindent
But (\ref{15}) and (\ref{16}) give $x^i_j=0$, against our assumption (\ref{13}). Now suppose that $2\leq j'<c_{i+k}+1$. In this case 
$$ x^i_jx^{i+k}_{j'+1}=0, $$
\noindent whence, by (\ref{13})
\begin{equation}\label{18} x^{i+k}_{j'+1}=0. \end{equation}
\noindent
Now, by (\ref{10}),
$$
M^{i+k}_{j'-1j'}({\bf x})=\left|\begin{array}{cc}
x^{i+k}_{j'-1}&x^{i+k}_{j'}\\
x^{i+k}_{j'}&x^{i+k}_{j'+1}
\end{array}
\right|
=x^{i+k}_{j'-1}x^{i+k}_{j'+1}-({x^{i+k}_{j'}})^2=0,
$$
\noindent
which, by (\ref{18}), yields $x^{i+k}_{j'}=0$. This, again, contradicts (\ref{13}). Hence (\ref{11}) is always true. This completes the proof.
\begin{example}\label{example3}{\rm  According to Theorem \ref{theorem2}, the ideal $J$ defined in Example \ref{Example2} has pure height 7 and is a s.t.c.i. on the following polynomials:
\begin{eqnarray*}
F^2_1&=&X_3X_5-X_4^2\\
F^3_1&=&X_5X_7-X_6^2\\
F^3_2&=&X_5X_8^2-2X_6X_7X_8+X_7^3\\
F^4_1&=&X_9X_{11}-X_{10}^2\\
G_3&=&X_1X_{11}\\
G_2&=&X_1X_8+X_3X_{11}\\
G_1&=&X_1X_5+X_3X_8+X_5X_{11}.
\end{eqnarray*}
Computations with  CoCoA \cite{cocoa} show that, for all fields $K$, the least power of $J$ contained in the ideal $(F^2_1,
F^3_1, F^3_2, F^4_1, G_1, G_2, G_3)$ is the 13th. 
}
\end{example}

\section{The fiber cone of a monomial variety of codimension two}
In this section we show that Theorem \ref{theorem2} applies to a relevant class of polynomial ideals. We first need to recall some preliminary notions from commutative algebra. Let $I$ be an ideal of $R$, and let $t$ be an indeterminate over $R$. The graded ring $R[It]=\oplus_{i\in{\bf N}}I^it^i$ is called the {\it Rees algebra} of $I$. If ${\cal M}$ is the ideal of $R$ generated by the set of indeterminates \underline{$X$}, so that $R/{\cal M}\simeq K$, then the quotient ring $F(I)=R[It]/MR[It]\simeq R[It]\otimes K$ is called the {\it fiber cone algebra}  of $I$. Suppose that $H_1,\dots, H_m$ form a minimal generating set of $I$, and introduce  a set of $m$ independent variables over $R$, say \underline{$T$}=$\{T_1,\dots, T_m\}$. Consider the ring homomorphism 
$$\phi:R[\mbox{\underline{$T$}}]\to R[It]$$
\noindent
such that, for all $i=1,\dots, m$,
$$\phi(T_i)=H_it.$$
\noindent 
It is evidently surjective, so that $R[It]\simeq R[\mbox{\underline{$T$}}]/\mbox{Ker}\,\phi$. Tensoring with $K$ yields $F(I)\simeq K[\mbox{\underline{$T$}}]/J$, for a suitable presentation ideal $J\subset K[\mbox{\underline{$T$}}]$.\par\smallskip\noindent
  Assume that $I\subset R=K[X_1,\dots, X_n, Y_1, Y_2]$ is the defining ideal of a {\it a  monomial variety of codimension two}, i.e., of a variety of $K^{n+2}$ admitting a parametrization of the following form:
$$x_1=u_1^{a_1},\ x_2=u_1^{a_2},\dots, \ x_n=u_1^{a_n}, \ y_1=u_1^{b_1}u_2^{b_2}\cdots u_n^{b_n}, \ y_2=u_1^{c_1}u_2^{c_2}\cdots u_n^{c_n},$$
\noindent
where, $a_1,a_2,\dots, a_n$ are positive integers and, for all $i=1,\dots, n$, the exponents $b_i, c_i$ are non negative integers such that $(b_i, c_i)\ne(0,0)$,  and, moreover, $(b_1,b_2,\dots, b_n)\ne(0,\dots, 0)$ and $(c_1,c_2,\dots, c_n)\ne(0,\dots, 0)$. This is an example of affine toric variety. If the above parametrization is homogeneous, it defines a projective variety of ${\bf P}^{n+1}$. 
The presentation ideal $J$ of the fiber cone algebra $F(I)$ is known to be of the type discussed in Section 1 (see \cite{BM1}, Proposition 3.6). An explicit construction of the barred matrix $A$ associated with $J$ can be found in \cite{G} or \cite{GMS}. We present an example which was considered, from a different point of view, in \cite{BM1}, Example 3.7 (b).  
\begin{example}\label{example4} {\rm Consider the projective monomial curve of ${\bf P}^3$ parametrized by
$$x_1=u_1^{534}, \ x_2=u_2^{534}, \ y_1=u_1^{245}u_2^{289}, \ y_2=u_1^{144}u_2^{390}$$
\noindent
Its defining ideal $I\subset R=K[X_1,X_2, Y_1, Y_2]$ is minimally generated by the following 6 binomials:
\begin{eqnarray*}
P_1&=&y_1^{42}-x_1^{19}x_2^{22}y_2\\
P_2&=&y_1^{12}x_2^{3}-x_1^2y_2^{13}\\
P_3&=&y_1^{30}y_2^{12}-x_1^{17}x_2^{25}\\
P_4&=&y_1^{18}y_2^{25}-x_1^{15}x_2^{28}\\
P_5&=&y_1^6y_2^{38}-x_1^{13}x_2^{31}\\
P_6&=&y_2^{51}-y_1^6x_1^{11}x_2^{34},
\end{eqnarray*}
\noindent
and the presentation ideal $J\subset K[T_1,\dots, T_6]$ of $F(I)$ is associated with the barred matrix
$$\left(\begin{array}{cccccc}
T_1&\vline\,\vline&T_3&T_4&\vline\,\vline&T_5\\
T_2&\vline\,\vline&T_4&T_5&\vline\,\vline&T_6
\end{array}
\right).
$$
\noindent
Therefore $J=(T_3T_5-T_4^2, \ T_1T_4,\ T_1T_5,\ T_1T_6,\ T_3T_6,\ T_4T_6,)$, and $J$ is a s.t.c.i. on the following 3 polynomials:
\begin{eqnarray*}F_1&=&T_3T_5-T_4^2\\
G_2&=&T_1T_6\\
G_1&=&T_1T_5+T_3T_6.
\end{eqnarray*}
It can be easily checked that $I^2\subset (F_1, G_1, G_2)$.}
\end{example}
\section{On varieties of minimal degree}
The classifcation of projective varieties of minimal degree (i.e., of degree equal to the codimension plus one, see \cite{EG} for the details) is due to the contributions of various authors. The irreducible ones are the quadric hypersurfaces and the cones over the Veronese surface in ${\bf P}^5$. The reducible case was settled by Xamb\'o, who proved the following
\begin{theorem}\label{theorem3}{\rm [\cite{X}, Section 1]}
Let $J$ be a reduced homogeneous ideal of $R$, defining a reducible variety $V$ of ${\bf P}^n$ of pure dimension $d$ and degree $\delta\leq n-d+1$. Then ${\bf P}^n$ contains linear subspaces $L_1,\dots, L_r$ and there are $d$-dimensional scrolls $V_i\subset L_i$ such that $V=\bigcup_{i=1}^rV_i$ and for each $i=2,\dots, r$, 
$$V_i\cap(V_1\cup\cdots\cup V_{i-1})=L_i\cap(\overline{L_1\cup\cdots\cup L_{i-1}}),$$
\noindent
which is a linear subspace of dimension $d-1$ (the bar denotes the projective closure). Moreover, $\delta=n-d+1$. 
\end{theorem}
The original version of the theorem contained the assumption of connectivity in codimension one, but this can be dropped, as shown in \cite{BM3}. A constructive characterization of the defining ideals of the varieties in Theorem \ref{theorem3} was given in \cite{BM2}. They include {\it all} the ideals studied in \cite{BM1}, i.e., not only the ideals of the type discussed in Section 1, but also the ideals associated to a more general kind of barred matrix,
\begin{eqnarray*}&&\!\!\!\!\!\!\!\!A=\\
&&\!\!\!\!\!\!\!\!\!\!\!\!\left(\begin{array}{ccccccccccccccccc}
B_{11}&\vline& B_{12}&\vline\ \cdots\ \vline& B_{1s_1}&\vline\,\vline&B_{21}&\vline\ \cdots\ \vline& B_{2s_2}&\vline\,\vline\ \cdots\ \vline\,\vline&B_{r1}&\vline\ \cdots\ \vline& B_{rs_r}
\end{array}\right),
\end{eqnarray*}
\noindent
 where every big block $B_i$ may consist of various {\it small blocks} $B_{i1}, \dots,B_{is_i}$, which are scroll matrices and have pairwise distinct sets of entries. We show that the s.t.c.i. property established in Theorem \ref{theorem2} extends to some, but not to all ideals associated with this larger class of barred matrices. We preliminarily remark that Theorem \ref{theorem1} is true also in this more general case (see \cite{BM1}, Corollary 1.3). 
\begin{example}\label{example4'}{\rm Let 
$$A=\left(\begin{array}{ccccc}
X_1&\vline& X_3&\vline\, \vline& X_4\\
X_2&\vline& X_4&\vline\, \vline& X_5
\end{array}\right).$$
\noindent
The associated ideal is 
$$J=(X_1X_4-X_2X_3, \ X_1X_5, \ X_3X_5),$$
\noindent
and height\,$J=2$. We show that $J$ is a s.t.c.i. on the following two polynomials
\begin{eqnarray*} P_1&=&(X_1X_4-X_2X_3)X_2+X_1X_5\\
P_2&=&(X_1X_4-X_2X_3)X_4+X_3X_5,
\end{eqnarray*}
\noindent
i.e., that
\begin{equation}\label{rad}J=\sqrt{(P_1,P_2)}.\end{equation}
\noindent
It holds:
\begin{eqnarray*}
(X_1X_5)^2&=&X_1(X_5-X_2X_4)P_1+X_1X_2^2P_2\\
(X_3X_5)^2&=&X_3(X_5+X_2X_4)P_2-X_3X_4^2P_1,
\end{eqnarray*}
\noindent
so that 
\begin{equation}\label{rad1}X_1X_5,\, X_3X_5\in \sqrt{(P_1,P_2)},\end{equation}
whence, by definition of $P_1$ and $P_2$,
$$(X_1X_4-X_2X_3)X_2, \ (X_1X_4-X_2X_3)X_4\in \sqrt{(P_1,P_2)}.$$
\noindent
It follows that 
$(X_1X_4-X_2X_3)^2\in\sqrt{(P_1,P_2)}$, i.e.,
\begin{equation}\label{rad2}X_1X_4-X_2X_3\in\sqrt{(P_1,P_2)}.\end{equation}
\noindent
Relations (\ref{rad1}) and (\ref{rad2}) imply (\ref{rad}).}
\end{example}
In the next example, the s.t.c.i. property fails to be true.
\begin{example}\label{example5}{\rm Let 
$$A=\left(\begin{array}{ccccccc}
X_1&\vline& X_3&\vline& X_5&\vline\,\vline&X_6\\
X_2&\vline& X_4&\vline& X_6&\vline\,\vline&X_7
\end{array}\right).$$
\noindent
The associated ideal is 
$$J=(X_1X_4-X_2X_3,\ X_3X_6-X_4X_5, \ X_1X_6-X_2X_5, \ X_1X_7, \ X_3X_7, \ X_5X_7),$$
\noindent
 and height\,$J=3$. 
}
\end{example} Let ara\,$J$ denote the arithmetical rank of $J$. 
We show that   ara\,$J>\,$height\,$J$.  We shall use the following criterion, which is based on \'etale cohomology and is due to Newstead \cite{N}.
\begin{lemma}\label{Newstead}{\rm [\cite{BS}, Lemma 3$^\prime$]} Let
$W\subset\tilde W$ be affine varieties. Let $d=\dim\tilde
W\setminus W$. If there are $s$ equations $F_1,\dots, F_s$ such
that $W=\tilde W\cap V(F_1,\dots,F_s)$, then 
$$\het^{d+i}(\tilde W\setminus W,{\bz}/r{\bz})=0\quad\mbox{ for all
}i\geq s$$ and for all $r\in{\bz}$ which are prime to {\rm char}\,$K$.
\end{lemma}
We refer to \cite{Mi} for the basic notions on \'etale cohomology. 
Let $p$ be a prime such that $p\ne$\,char\,$K$. Let $V=V(J)\subset K^7$. In view of Lemma \ref{Newstead}, for our purpose it suffices to show that 
\begin{equation}\label{19} \het^{10}(K^7\setminus V, {\bz}/p{\bz})\ne0.
\end{equation}
\noindent
 By Poincar\'e Duality (see \cite{Mi}, Theorem 14.7, p.~83) we have
\begin{equation}\label{20} {\rm Hom}_{\scriptstyle{\bz}/p{\bz}} (\het^{10}(K^7\setminus V, {\bz}/p{\bz}), {\bz}/p{\bz})\simeq  \hc^{4}(K^7\setminus V, {\bz}/p{\bz}),
\end{equation}
\noindent
where $\hc$ denotes \'etale cohomology with compact support. 
For the sake of simplicity, we shall omit the coefficient group ${\bz}/p{\bz}$ henceforth. Let $W$ be the subvariety of $K^7$ defined by 
$$X_1X_4-X_2X_3=0,\ X_3X_6-X_4X_5=0,\ X_1X_6-X_2X_5=0, \ X_7=0.$$
\noindent
Then $W\subset V$, and 
$$V\setminus W=\{(x_1,\dots, x_7)\in K^7\vert x_1=x_3=x_5=0, \ x_7\ne0\}\simeq K^3\times (K\setminus\{0\}).$$
\noindent
It is well-known that
\begin{equation}\label{21}\hc^i(K^t)\simeq\left\{\begin{array}{cl} {\bz}/p{\bz}&\mbox{if }i=2t\\
0&\mbox{else, }
\end{array}\right.
\end{equation}
\noindent
and 
\begin{equation}\label{22}\hc^i(K^t\setminus\{0\})\simeq\left\{\begin{array}{cl} {\bz}/p{\bz}&\mbox{if }i=1,2t\\
0&\mbox{else. }
\end{array}\right.
\end{equation}
\noindent
Moreover, by the K\"unneth formula (\cite{Mi}, Theorem 22.1), 
$$\hc^i(V\setminus W)\simeq\displaystyle\bigoplus_{h+k=i}\hc^h(K^3)\otimes\hc^k(K\setminus\{0\}),$$
\noindent
so that, by (\ref{21}),
\begin{equation}\label{23}\hc^3(V\setminus W)=\hc^4(V\setminus W)=0.\end{equation}
We have a long exact sequence of \'etale cohomology with compact support:
$$\cdots\rightarrow\hc^3(V\setminus W)\rightarrow \hc^4(K^7\setminus V)\rightarrow \hc^4(K^7\setminus W)\rightarrow 
\hc^4(V\setminus W)\rightarrow\cdots.$$
\noindent
By (\ref{23}) it follows that
\begin{equation}\label{24} \hc^4(K^7\setminus V)\simeq\hc^4(K^7\setminus W)
\end{equation}
\noindent
In view of (\ref{20}) and (\ref{24}) our claim (\ref{19}) will follow once we have proven that 
\begin{equation}\label{25} \hc^4(K^7\setminus W)\ne0.\end{equation}
\noindent
 This is what we are going to prove next. Note that $W$ is the variety of $K^6$ defined by the vanishing of the 2-minors of a generic $2\times 3$ matrix of indeterminates. Thus $\bar W=W\setminus\{0\}$ is the set of non zero 2$\times$3 matrices with proportional rows. The set of such matrices where the first row is zero is a closed subset of $\bar W$ which can be identified with $Z=K^3$, and its complementary set is $\bar W\setminus Z\simeq K^3\setminus\{0\}\times K$.  We thus have a
long exact sequence of \'etale cohomology with compact support:
\begin{equation}\label{27a}\cdots\rightarrow\hc^2(Z)\rightarrow \hc^3(\bar W\setminus Z)\rightarrow \hc^3(\bar W)\rightarrow \hc^3(Z)\rightarrow\cdots,\end{equation}
\noindent
where, according to (\ref{21}), $\hc^2(Z)=\hc^3(Z)=0$, whereas, by  the K\"unneth formula and (\ref{22}), 
$$\hc^3(\bar W\setminus Z)\simeq \hc^1(K^3\setminus\{0\})\otimes\hc^2(K)\simeq \bz/p\bz.$$
\noindent
It follows that (\ref{27a}) gives rise to an isomorphism
\begin{equation}\label{26}\hc^3(\bar W)\simeq{\bz}/p{\bz}.\end{equation}
\noindent
On the other hand there is also the following long exact sequence of \'etale cohomology with compact support:
\begin{equation}\label{27}\cdots\rightarrow\hc^3(K^7\setminus \{0\})\rightarrow \hc^3(\bar W)\rightarrow \hc^4(K^7\setminus W)\rightarrow 
\cdots,\end{equation}
\noindent
where 
\begin{equation}\label{28} \hc^3(K^7\setminus\{0\})=0.\end{equation}
\noindent
by (\ref{22}). Hence, in view of (\ref{28}) and (\ref{26}), in (\ref{27}) the left term is zero, and the middle term is non zero. It follows that the right term is non zero, i.e., claim (\ref{25}) holds. This proves that ara\,$J>3$, as desired.
It can be easily checked that the variety $V$ is defined by the following 5 equations:
$$X_1X_4-X_2X_3=0,\ X_3X_6-X_4X_5=0,$$
$$X_1X_6-X_2X_5+X_3X_7=0, \ X_1X_7=0, \ X_5X_7=0.$$ 
\noindent
Hence ara\,$J\leq5$. We conjecture that ara\,$J=5$. 
\par\medskip\noindent
\'Etale cohomology is not the only possible tool for finding a lower bound for the arithmetical rank. In general we have
\begin{equation}\label{29}{\rm cd}\,J\leq{\rm ara}\,J,\end{equation}
\noindent where 
$${\rm cd}\,J=\max\{i\vert H^i_J(R)\ne0\},$$
\noindent 
is called the {\it cohomological dimension} of $J$. Here $H^i_J$ denotes the $i$-th local cohomology group with respect to $J$;  we refer to Huneke \cite{Hu} for an extensive exposition of this subject. Inequality (\ref{29}) can be strict, but not many examples  of this kind are known so far. One monomial ideal was found by Zhao Yan \cite{Z}, Example 2 following \cite{SV}, p.~250, and \cite{L}, Example 1,  whereas classes of determinantal ideals are described in \cite{BS} and \cite{B0}. In all these cases (\ref{29}) is strict in all but one characteristics. We show that for the ideal $J$ discussed in Example \ref{example5} we always have cd\,$J<5$: if our above conjecture is true, this would provide an example of ideal whose cohomological dimension differs from its arithmetical rank in all characteristics.  According to \cite{Hu}, Theorem 2.2, we have the following long exact sequence of local cohomology:
\begin{equation}\label{30} 
\cdots\rightarrow H^i_{J+(X_7)}(R)\rightarrow H^i_J(R)\rightarrow H^i_J(R_{X_7})\rightarrow\cdots,\end{equation}
\noindent
where, by \cite{Hu}, Proposition 1.10, $H^i_J(R_{X_7})\simeq H^i_{J_{X_7}}(R_{X_7})$. Now
\begin{eqnarray*} J+(X_7)&=&(X_1X_4-X_2X_3,\ X_1X_6-X_2X_5,\ X_3X_6-X_4X_5,\ X_7)\subset R,\\
J_{X_7}&=&(X_1,\ X_3,\ X_5)\subset R_{X_7}.
\end{eqnarray*}
\noindent
Since $J+(X_7)$ and $J_{X_7}$ are generated by sets of at most 4 elements in $R$ and $R_{X_7}$ respectively, inequality (\ref{29}) implies that
$$H^i_{J+(X_7)}(R)=H^i_{J_{X_7}}(R_{X_7})=0,\qquad\mbox{for all }i\geq5,$$
\noindent
so that, in view of (\ref{30}), 
$$H^i_J(R)=0\qquad\mbox{for all }i\geq5.$$
\noindent
Hence cd\,$J<5$, as was to be shown.  
\begin{remark}{\rm According to \cite{EG}, Theorem 4.2, every variety fulfilling the assumption of Theorem \ref{theorem3} is Cohen-Macaulay. By \cite{PS}, Prop. (4.1), this implies that cd\,$J=\,$ht\,$J=3$. Hence we certainly have that cd\,$J<$ara\,$J$ in all positive characteristics. 
}\end{remark}

\end{document}